\documentclass[12pt]{article}

\usepackage[utf8]{inputenc}
\usepackage{graphicx}
\usepackage{color}
\usepackage[english]{babel}
\usepackage{here}
\usepackage{amsmath, amsthm, amsfonts}

\newtheorem{thm}{Theorem}[section]

\theoremstyle{definition}

\theoremstyle{remark}
\newtheorem{rem}[thm]{Remark}

 \usepackage{vmargin}

 \setmargins
 {1.5cm}       % margen izquierdo
 {0cm}                        % margen superior
 {18cm}                      % anchura del texto
 {24cm}                    % altura del texto
 {10pt}                           % altura de los encabezados
 {1cm}                           % espacio entre el texto and los encabezados
 {0pt}                             % altura del pie de pogina
 {2cm}                           % espacio entre el texto and el pie de pogina
 
 \title{Finite difference and finite element methods for partial differential equations on fractals}
 \author{Luis F. Contreras H. \and Juan Galvis
 	\thanks{{\tt lfcontrerash,jcgalvisa@unal.edu.co}, Departamento de Matem\'aticas, Universidad Nacional de Colombia, Carrera 45 No. 26-85, Edificio Uriel Gutierr\'ez, Bogot\'a D.C., Colombia.}
 	\date{}
 }
 
\begin{document}

\maketitle

\begin{abstract}
In this paper, we present numerical procedures to compute solutions of partial 
differential equations posed on fractals. In particular, we consider the strong form of the equation using standard graph Laplacian matrices and also weak forms of the equation derived using standard length or area measure on a discrete approximation of the fractal set.   We then introduce a numerical procedure to normalize the obtained diffusions, that is, a way to compute the renormalization constant needed in the definitions of the actual partial differential equation on the fractal set. A particular case that is studied in detail is the solution of the Dirichlet problem in the Sierpinski triangle. Other examples are also presented including a non-planar Hata tree.

\end{abstract}

\section{Introduction}

In recent years we have seen many applications of fractal sets in modeling sciences. Especially, to study several processes that can be modeled using fractals and specially self-similar structures. We mention processes related to diffusion on fractal sets witch have several possible applications, including diffusion of substances in biological structures and flow inside fractures in modeling fluid flow in fractured porous media, among other models. See \cite{JK,RS1}. In this paper, we consider fractals defined as self-similar sets with the additional properties  of being post-critically finite;
see \cite{JK}. These self-similar sets can be approximated (in the Hausdorff metric) by a finite union of sets generated by removing 
a finite number of vertices from a graph approximation of the fractal set.

We recall the definition of the Laplace operator, a standard model for fractal diffusion. Then, introduce numerical approximation procedures of the presented model. It is important to stress that our approximation procedure consists of renormalizing standard approximation methods 
on two- and three-dimensions such as the finite difference and the finite element method. 
See for example \cite{GRS} where a finite element method is designed and analyzed. Note that in \cite{GRS}, the authors assume the renormalization constant is known. We do not need to compute the renormalization constant analytically and instead of of that, we approximate it numerically.
We also mention that the finite element formulation implemented here is based on weak forms computed using standard length and area measures restricted to approximations of the fractal sets. In particular, in the case of the 
Sierpinski triangle, the approximation by the finite element method can be summarized as follows: 
\begin{itemize}
\item The computations are carried out on approximations of the fractal (that could be a union of edges or triangles).
\item Approximation of the computations of derivatives for which we use classical derivatives of piecewise linear functions in one and two dimensions. Alternatively, we also approximate derivatives using standard weight and adjacency graph matrices.
\item Approximation of the self-similar measure. Here we test different approximations: 
1) The measure induced by the length measure  restricted to the edges of the triangles in the finite graph that represents the current approximation of the fractal, and 2) The measure induced by the area measure restricted to the triangles of the current approximation of the fractal.
\item Approximation of the values for the rescaling or renormalizing to obtain renormalized operators and guarantee that the obtained solution approximates the solution of the continuous problem. 

\end{itemize}

In the last step, the finite element procedure above is then renormalized with a pre-computation of the renormalization constant to obtain the correct approximation of the Laplace operator in the Sierpinski triangle.
The computation of the renormalization constant involves the comparison of two (not-normalized) solutions at consecutive (levels of refinement) approximations of the fractal set. 
We call this procedure the renormalized FEM, or rFEM for short. We illustrate the performance of the rFEM numerically for the solution of a Dirichlet problem on the Sierpinski triangle and a realization of the Hata tree. We also present a renormalized Finite Difference  (rFD) procedure using the same idea.

The rest of the paper is organized as follows. 
In Section \ref{sec:selsimilar} we recall
some examples of self-similar sets. 
Section \ref{sec:graphlap}
is dedicated to reviewing the definition of the graph laplacian operator. 
In Section \ref{sec:dirichlet} we present several formulations of the Dirichlet problem as well as the proposed procedure for the computation of the renormalization constant. We finish this section with numerical experiments and illustrations of the correctness of our method.  
In Section \ref{sec:conclusions} we present some conclusions. 

\section{Examples of self-similar sets}\label{sec:selsimilar}
This section reviews some facts related to self-similar sets, its constructions and also related to the approximation of fractal sets. We follow 
\cite{MJ, JK,RS1}.

A self similar set is obtained by applying  a fixed point functional iteration. Let  $(X,d)$ be a Hausdorff metric space and denote by 
$C(X)$ the metric space  of all compact subsets of a metric space equipped in the Hausdorff metric.  Assume  you have contractions $f_i:X\rightarrow X$, $i=1,2,3,...,N$ and define $F:C(X)\rightarrow C(X)$ by $F(A)=\bigcup_{1\leq i\leq N} f_i(A)$ for all $A\in C(X)$. Then $F$ have a unique fixed point $K$. Also for any $A\in C(X)$, $F^n(A)$ converge to $K$ when $n\rightarrow \infty$ with respect to Hausdorff metric. We then have that
there  exists a unique non-empty compact $K\subset X$ such that
\begin{equation}\label{eq:pfF}
K=f_1(K)\cup f_2(K)\cup...\cup f_N(K).
\end{equation}
The set $K$ is the self-similar set
associated with $\left\{f_1,f_2,\ldots,f_N\right\}$ and this is a fractal set.
The term self-similar is given to $K$ 
because $K$ is the union of images of itself by the contractions. See \cite{JK} and reference therein.

In order to obtain particular examples, 
we need only the initial metric and 
the finite set of contractions. In 
particular, we can consider the following 
examples of subsets of  $\mathbb{R}^d$, $d=1,2,3$ with the Euclidean distance together with a finite family of contractions. See  \cite{MJ, JK,RS1}.

\begin{itemize}
\item {The Kosh curve:} Let  $a_1=(0,0)$ and $a_2=(1,0)$ be the initial nodes for the construction. This set is denoted with $W_0=[a_1,a_2]$\footnote{For $a,b\in \mathbb{R}^2$, $[a,b]$ denote the line segment from $a$ to $b$.}. Consider the 
contractions: 
$$f_i(x)=\frac{1}{3}r(\theta)^i\cdot x+\left(\frac{i}{3},0\right) \quad \mbox{ where }\quad i=0,1,2,3.$$
Here  $r(\theta)$ is the rotation matrix with angle $\theta$. We define the set $\left\{W_n: n\in \mathbb{N}_0\right\} \subseteq \mathbb{R}^2$ inductively by,
$$
W_{n+1}=f_0(W_n)\cup f_1(W_n)\cup f_2(W_n)\cup f_3(W_n).
$$
We can define $K=\lim_{n\rightarrow \infty} W_n$ with the limit in the Hausdorff metric. The set $K$ is called the Kosh curve. In particular, $W_n$ is an approximation of the Kosh curve where discrete differential operators can be constructed in order to 
approximation differential operators defined 
on the Kosh curve $K$.

\begin{figure}[H]\label{fig:Koch}
	\begin{center}
		\includegraphics[scale=0.4]{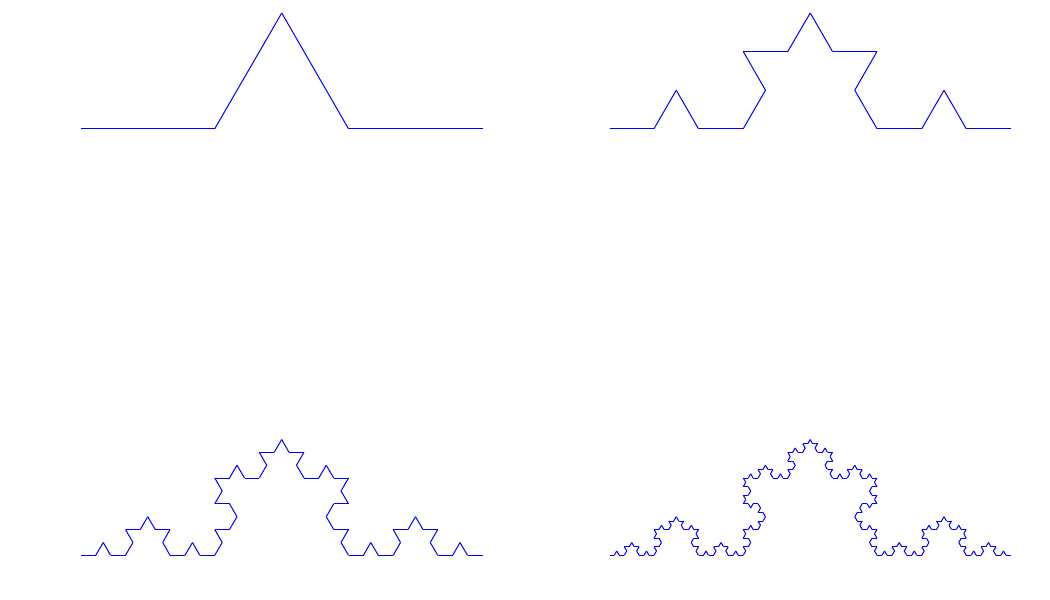}
		\caption{{\small The first four iterations of the construction off the Kosh curve: $W_1$ (Up right), $W_2$ (Up left), $W_3$ (down right), $W_4$ (down left).}}
	\end{center}
\end{figure}

\item{The Sierpinski triangle:} 
The Sierpinski triangle is one of the most known examples of self-similar, see \cite{CJ}. 
It can be consider a benchmark fractal where several questions and problems can be test out. A construction of the triangle goes as 
follows. Let  $a_0=(0,0)$, $a_1=(1,0)$, $a_2=(\frac{1}{2},\frac{\sqrt{3}}{2})$ the vertices of the equilateral triangle  $X \subseteq \mathbb{R}^2$. We consider the set $X$ with the Euclidean distance. For each  $i=0,1,2$  we define the affine mapping
\begin{eqnarray}
f_i: X &\rightarrow & X,\\
x &\rightarrow & f_i(x)=\frac{1}{2} (x-a_i)+a_i.
\end{eqnarray}
Let $W_0=[a_0,a_1]\cup [a_1,a_2]\cup[a_2,a_0]$ and we define   $\left\{ W_n: n \in \mathbb{N}_0 \right\}$ by
$$W_{n+1}=f_0( W_n)\cup f_1( W_n)\cup f_2( W_n) \mbox{ for all } n\in \mathbb{N}_0.$$

As before, we see that $ K=\lim_{n\to \infty} W_n$ (where the limit is taken in the Hausdorff metric), and $W_n$ can be viewed as an approximation of $K$ where differential operators can be computed to approximation differential operators defined on $K$.

\begin{figure}[h]
	\begin{center}
		\includegraphics[scale=0.6]{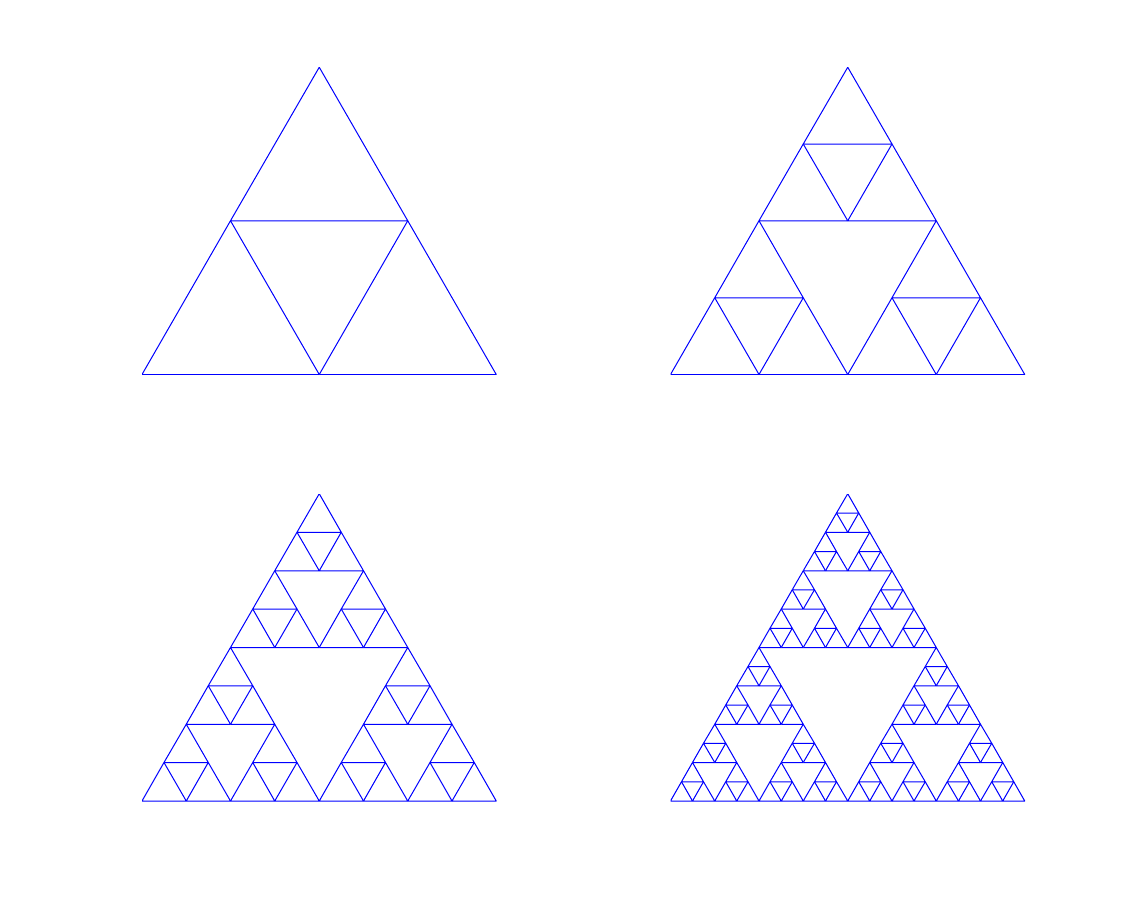} 
		\caption{{\small The first four iterations of the Sierpinski triangle: $W_1$ (up right), $W_2$ (up left), $W_3$ (up right), $W_4$ (up left).}}\label{k4}
	\end{center}
\end{figure}

\item{The Hata tree in the plane: }
Let $p_1=(0,0)$ and $p_2=(1,0)$ the vertices of $W_0=[p_1,p_2]$. Define
\begin{eqnarray*}
	f_0(x)&=&\frac{x}{3},\\
	f_1(x)&=&\left( \frac{1}{3},0\right)+\frac{x}{3}\cdot r(\pi/3),\\
	f_2(x)&=&\left( \frac{1}{3},0\right)+\frac{x}{3},\\
	f_3(x)&=&\left( \frac{2}{3},0\right)+\frac{x}{3}\cdot r(\pi/3),\\
	%\mbox{y }\\
	f_4(x)&=&\left( \frac{2}{3},0\right)+\frac{x}{3}.
\end{eqnarray*}

Where $r(\theta)$ denotes the $\theta$-rotation matrix. Define  
$$W_{n+1}= \bigcup_{i=0}^4 f_{i}( W_n) \mbox{ for all } n.$$

The first four iterations are illustrated in Figure  \ref{HTree}.

\begin{figure}[H]
	\begin{center}
		\includegraphics[scale=0.6]{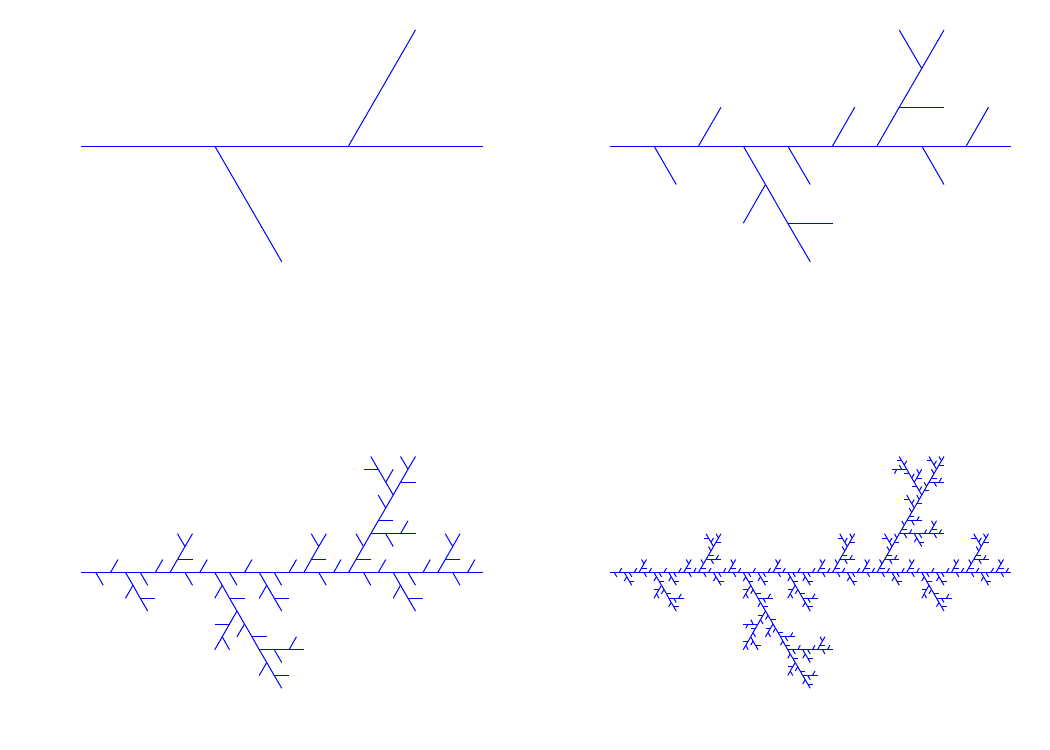} 
		\caption{{\small The first four iteration of the Hata tree $W_1$ (up right), $W_2$ (up left), $W_3$ (down right) and $W_4$ (down left).}}\label{HTree}
	\end{center}
\end{figure}

\item  {Non-planar Hata tree:} Let $p_1=(0,0,0)$ and $p_2=(0,0,1)$ the vertices of the set $W_0=[p_1,p_2]$. Define,
\begin{eqnarray*}
	f_0(x)&=&\frac{x}{3},\\
	f_1(x)&=&\left(0,0, \frac{1}{3}\right)+k_2\frac{x}{3}+\frac{k_1\cdot n_1}{2}+\frac{\sqrt{3}}{2}\cdot k_1\cdot n_2,\\
	f_2(x)&=&\left(0,0 ,\frac{1}{3}\right)+\frac{x}{3},\\
	f_3(x)&=&\left(0,0, \frac{2}{3}\right)+k_2\frac{x}{3}-\frac{1}{2}\cdot k_1\cdot n_1+\frac{\sqrt{3}}{2}\cdot k_1\cdot n_2,\\
	f_4(x)&=&\left(0,0, \frac{2}{3}\right)+k_2\frac{x}{3}-k_1\cdot n_1,\\
	f_5(x)&=&\left(0,0, \frac{2}{3}\right)+\frac{x}{3}.\\
\end{eqnarray*}

Where $k_1=\frac{1}{3}\sin(\pi/4)$, $k_2=\frac{1}{3}\cos(\pi/4)$ and $\left\{n_1, n_2,\frac{x}{|x|}\right\}$  is and orthonormal set.

\begin{figure}[H]
	\begin{center}\label{HTree3}
		\includegraphics[scale=.7]{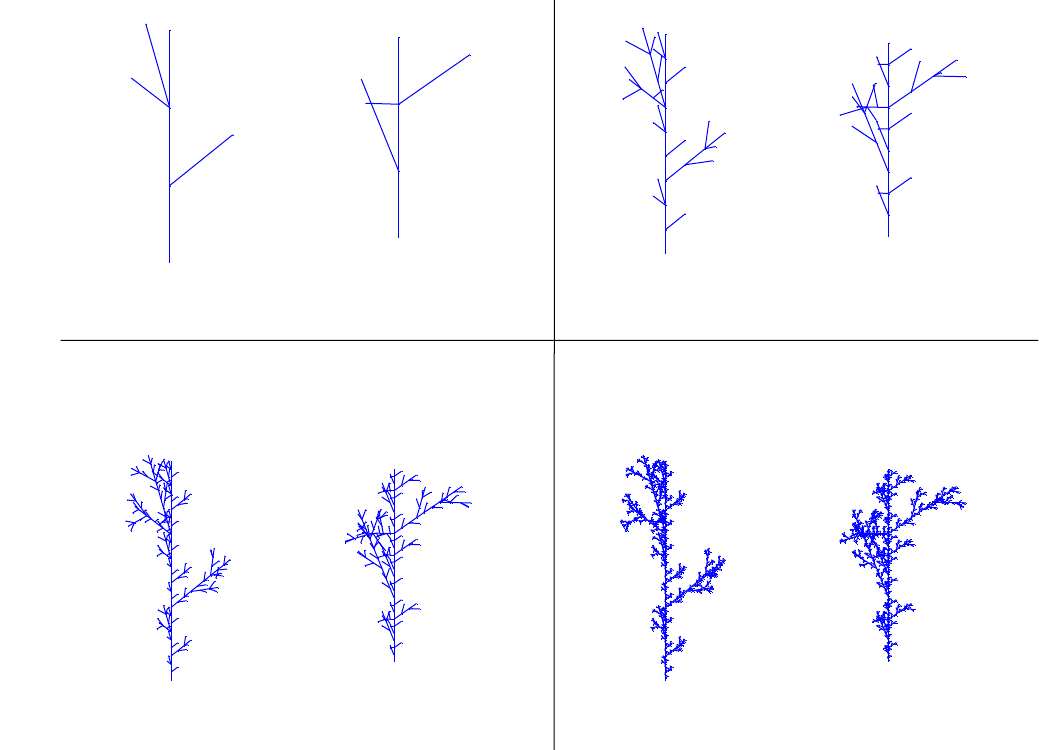} 
		\caption{{\small The firstfour iteration of the Hata tree. See  $W_1$ (Up right), $W_2$ (Up left), $W_3$ (Down right) and $W_4$ (Down Left).
			}}
		\end{center}
	\end{figure}
\end{itemize}

The previous fractals can be built from the set of initial vertices $V_0$ and use the sequence $V_{n+1}=f_1(V_n)\cup...\cup f_N(V_n)$. We will use this notation later on. For example, for the Sierpinski triangle the initial vertices are $V_0=\left\{ a_1,a_2,a_3\right\}$ and we will have the sequence $V_{n+1}=f_1(V_n)\cup...\cup f_3(V_n)$.

\section{ Laplacian  on  a  graph}\label{sec:graphlap}
In this section, we review the construction of the Laplace operator and the energy on a graph; in particular, we introduce the renormalization constant, which is important to get a  finite limit of the energies associated with a family of graphs that approximates a fractal.  
See \cite{ MJ, JK,RS1}.

Let $G(V,E)$ a finite graph, where  $V=\{v_1,v_2,...,v_n\}$ determines the set of vertices and $E$ the set of edge (without orientation) of $V$. If  $v,w\in V$ and exist an edge between  $v$ and $w$ we write $v\sim w\in E$. define the adjacency matrix $A_G$ associate to a graph $G$  as the $n\times n$  matrix 
	$ A_G=[a_{ij}]_{i,j=1}^n$, where
	\[
	a_{ij}= \left\{ \begin{array}{lcl}
	1, & \mbox{ if }\quad  v_i\sim v_j\in E, \\
	&  \\
	0, & \mbox{ in other case. } 
	\end{array}
	\right.
	\]
The weight matrix $P_G$ of  $G$ is the diagonal matrix of dimension  $n\times n$ defined by $P_G=[p_{ij}]_{i,j=1}^n$ with
	$p_{ij}=0$ when  $i\not = j$ and  
	$p_{ii}$ is the number of adjacent vertices  to $v_i$, $i=1,\dots, n$. Therefore, 
	\begin{equation}
	p_{ii} = \# \{ w \quad : \quad v\sim w\in E\}. 	\label{pesos} 
	\end{equation}
The Laplacian matrix  associate to    $G(V,E)$ is given by
	$$\Delta_G=P_G-A_G.$$
	\label{laplan}
If $u:V\to \mathbb{R}$ we define the energy of  $u$ by $$\mathcal{E}_G(u,u)=\sum_{x\sim y} (u(x)-u(y))^2.$$
The bilinear form associated to the energy 
is
$$\mathcal{E}_G(u,v)=\sum_{x\sim y}(u(x)-u(y))(v(x)-v(y)).$$

We denote $\mathcal{E}(u )=\mathcal{E}(u,u)$.
If we introduce the vectors  $\vec{U},\vec{V}\in \mathbb{R}^n$ given by  
$\vec{U}=\{ u(x)\}_{x\in V}$  and $\vec{V}=\{ v(x)\}_{x\in V}$, then 
\[
\mathcal{E}_G(u,v)= \vec{U}^T \Delta_G \vec{V} \quad \mbox{ and } \quad 
\mathcal{E}_G(u,u)= \vec{U}^T \Delta_G \vec{U}.
\]
We observe that matrix  $\Delta_G$ is the matrix representation of the energy $\mathcal{E}_G$. Now we consider the approximations $\{V_n\}_{n=1}^\infty$ to a fractal set $K$. We denote the energy associated to $V_n$ by 
\begin{equation}\label{defen}
\widetilde{\mathcal{E}}_n(u,v) = \mathcal{E}_{G_n}(U_n,V_n).
\end{equation}

The renormalized energy of level  $n=2,3,\dots$ can be computed for $u_n:V_n\to \mathbb{R}$ as
\begin{equation}\label{deferenor}
\mathcal{E}_n( u_n,u_n) = r^n \widetilde{\mathcal{E}}_n( u_n, u_n). 
\end{equation}
Here $r^n$ is a renormalization constant needed in order to obtain  a non-increasing sequence of renormalized energies $\mathcal{E}_n$, $n=0,1,\dots$. This step is necessary to obtain well-defined energy defined on $K$ that we introduce as a limit of the renormalized energies above. For more details  see \cite{JK} and references therein.

\subsection{The case of the Sierpinski triangle $K$}

For the Sierpinski triangle, see \cite{GRS}, we obtain that the energy for $f:V_\infty=\bigcup_{n} V_n \rightarrow \mathbb{R}$ is defined for $n\in \mathbb{N}_0$  by

\begin{equation}\label{rforEn}
\mathcal{E}_n(u,u)=\left(\frac{5}{3}\right)^n\sum_{x\sim_n y}(u(x)-u(y))^2,
\end{equation}
where $x\sim_n y$ is contained in $W_n$.
That is, the renormalization constant is $r=(5/3)$. Then, we can compute  $\mathcal{E}(u,u):=\lim_n \mathcal{E}_n(u|_{V_n},u|_{V_n})$. 
Introduce the renormalized  Laplace operators by \begin{equation}\label{deflaplacerenor}
\Delta_n= {5}^n  \widetilde{\Delta}_n.
\end{equation}

If we define the measure $\mu^n$ on $V_n$ by assigning full measure 1 to $V_n$ and stating that each cell has the same measure ($3^{-n}$), we see that we have $$\mathcal{E}_n(u,u)=-\int u \Delta_n ud\mu^n.$$

Using this identity it is defined the 
Laplace operator in $K$ by 
\begin{equation}
-\int u(x)\Delta u(x)\mu(dx)=\mathcal{E}(u,u). \label{deflaplacefrac} 
\end{equation}

Here $\mu$ is the standard self-similar measure associated to $K$ that can computed as the limit of the measure $\nu^n$ in the sense that

$$\int_K gd\mu =\lim_{n\to\infty}\int_{V_n} gd\nu^n,$$ 

where we note that, 
\begin{equation}
\int_{V_n} g d \nu^n=3^{-n} \left( \frac{2}{3} \sum_{x\in V_n\setminus V_0} g(x) + \frac{1}{3} \sum_{x\in V_0}g(x)\right).
\label{intauto}
\end{equation}

From here on, we will focus our study on the Sierpinski triangle.

\section{Formulations for the Dirichlet problem} \label{sec:dirichlet}

Given $g:V_{\infty}\rightarrow \mathbb{R}$, we seek for $u:V_\infty\rightarrow \mathbb{R}$ such that \begin{equation}
\begin{cases}
-\Delta u(x)=g(x) \quad x\in V_\infty\setminus V_0,\\
u(x)=h(x), 	\quad x  \in V_0, 
\label{fuerte}
\end{cases}
\end{equation}
where $\Delta u(x)$ is defined in \eqref{deflaplacefrac}. We call this the strong formulation of the Dirichlet problem.

Using the integration by parts formula we have $$-\int\Delta u\cdot vd\mu =\mathcal{E} (u,v). $$
Therefore, we can  write the problem as seeking for $u$ with bounded energy, 
$u\in H^1$, such that  
\begin{equation} \label{weakformJG}
\begin{cases} 
\mathcal{E} (u,v)=\displaystyle \int g vd\mu \quad \mbox{ for all }v\in H_0^1 \\
u(x)=h(x), \quad x  \in V_0.
\end{cases}
\end{equation}
Here, $H^1=\{u:V_\infty\rightarrow \mathbb{R},\quad \mathcal{E}(v,v)<+\infty \}$  and $H_0^1=\{ v\in H^1: v(x)=0,  x\in V_0 \}$. We refer to these formulation as 
the weak form of the Dirichlet problem.

\subsection{Finite difference approximation} \label{diferenciafinitas}

To approximate the solution of  (\ref{fuerte}) 
in $V_n$ we consider 
$\Delta_n $, the renormalized Laplace operator (defined for the Sierpinski triangle in \eqref{deflaplacerenor}).
Let the approximation be defined by $u^{FD}:V_n\to\mathbb{R}$
that can be written as 
$u^{\text{FD}}=\{ u^{FD}_n(x)\}_{x\in V_n}$ 
and we partitioned it as follows
\[
u^{\text{FD}}=\{ 
\{ u^{FD}_n(x)\}_{x\in V_0}, 
\{ u^{FD}_n(x)\}_{x\in V_n\setminus V_0} 
\} 
=[u_0,u_I].
\]
Note that $u_0$ is know and corresponds to the boundary values. Analogously for $b^{\text{FD}}=\{ g(x)\}_{x\in V_n} =
[b_0,b_I].$
We obtain the  block structure
\[
\Delta_n= \left(
\begin{array}{cc}
\Delta_{0,0} & \Delta_{0,I}\\
\Delta_{I,0} & \Delta_{I,I}\\
\end{array}\right).
\]
We compute $u_I$ as the solution of
\[
\Delta_{I,I} u_I = b_I - \Delta_{I,0}u_0.
\]

%Then  $U_0=\{ h(x)\}_{x\in V_0}$ and $U_I$ solve\begin{equation}\label{df}
%\Delta_{II} U_I =  g_I - \Delta_{I,0} U_0.
%\end{equation}
The main issue with this discrete formulation 
is that we need to know the renormalization constant:
the value $5^n$ in the case of the 
Sierpinski triangle, see \eqref{deflaplacerenor}. The renormalization constant can be viewed as a scaling of the forcing term in the linear system. This scaling can be approximated as explained next. 

\subsection{Computation of the renormalization constant}\label{procedure}
First we consider the finite difference method renormalization constant. The idea is to use the problem,
\begin{equation}
\begin{cases}
-\Delta u(x)=1 \quad x\in K,\\
u(x)=0,	\quad x  \in V_0. 
\label{consalte}
\end{cases}
\end{equation}
Denote by  $q^n$ and approximation of the 
renormalization factor that we want to compute. The space of study is $V_n/V_0$. We can approximate  (\ref{consalte}) by
\begin{equation*}
\begin{cases}
q^{n}\widetilde{\Delta}_n u_n(x)=1 \quad x\in V_{n}\setminus V_0,\\
u_n(x)=0,	\quad x  \in V_0.
\end{cases}
\end{equation*}
We want to compute the value $q^n$.
Numerically we can compute the solution of the problem, see \cite{SL},
\begin{equation*}
\begin{cases}
\widetilde{\Delta}z_n=1 \quad x\in V_{n}\setminus V_0,\\
z_n(x)=0, 	\quad x  \in V_0.
\end{cases}
\end{equation*} 
We must have that $z_n(x)=q^{n}u_n(x)$, $x\in V_n$ since $\widetilde{\Delta}_n$ is nonsingular. For $V_{n+1}$ we will have $z_{n+1} = q^{n+1}u_{n+1}$. Note that $z_n$ and $z_{n+1}$ can be calculated without knowing $q^{n}$. For  $n$ large enough we should have,
\begin{eqnarray}
z_n&=&q^{n}u_n \quad \mbox{ with }\quad u_n\approx u,\\
z_{n+1}&=&q^{n+1}u_{n+1} \quad \mbox{ with }\quad u_{n+1}\approx u,
\end{eqnarray}
where $u$ is the exact solution of \eqref{consalte}. Therefore we should be able to use the approximation,
\begin{equation}\label{definiqn}
q\approx q_{n,n+1}:=\frac{q^{n+1}u_{n+1}(x)}{q^{n}u_n(x)}=
\frac{z_{n+1}(x)}{z_{n}(x)},\end{equation}
where $x\in V_n\setminus V_0$. See some numerical illustration in  
Table \ref{tablaqn} for the case of the 
Sierpinski triangle. In this case we see hat $\lim_{n\to \infty} q_{n,n+1}=5$.\\

\begin{table}[h]
	\centering
	\begin{tabular}{ | c | c | c |}
		\hline
		$(n,n+1)$ & $\max_{x}q_{n,n+1}$ &
		$\mbox{mean}_xq_{n,n+1} $\\
		\hline
		(3,4)&  5  & 5\\
		(4,5)  &  5  &  5\\
		(5,6) &   5 &   5\\
		\hline
	\end{tabular}
	\caption{Values of  $q_{n,n+1}$ for the Sierpinski triangle, given for  \eqref{definiqn}.}
	\label{tablaqn}
\end{table}

A similar procedure can be implemented for the energy renormalization constants. The weak form of (\ref{consalte}) is written as, 
\begin{equation*}
\begin{cases}
r^{-n}\widetilde{\mathcal{E}}(u,v)=\int v d\mu_n \quad \mbox{ for all test } v,\\
u_n(x)=0 	\quad x  \in V_0.
\end{cases}
\end{equation*}

Compute the solution of the graph energy problem (without renormalization),
\begin{equation*}
\begin{cases}
\widetilde{\mathcal{E}}(z_n,v)= \int v d\mu_n \quad \mbox{ for all } v,\\
z_n(x)=0 	\quad x  \in V_0.
\end{cases}
\end{equation*} 

We can then define the approximation, \begin{equation}\label{definicionrn}
r=\frac{r^{-n}u_n(x)}{r^{-n-1}u_{n+1}(x)}\approx r_{n,n+1}= \frac{z_{n}(x)}{z_{n+1}(x)}.
\end{equation}
Here $x\in V_n\setminus V_0$. See a numerical 
verification in Table \ref{tablarn}. 
In this case we should have 
$\lim_{n\to \infty} r_{n,n+1}\approx 1.666...=\frac{5}{3}.$

\begin{table}
	\centering
	\begin{tabular}{ | c | c | c |}
		\hline
		$(n,n+1)$ & $\max_{x}r_{n,n+1}$ &
		$\mbox{mean}_xr_{n,n+1}$\\
		\hline
		(3,4) &  1.6667 & 1.6667 \\
		(4,5) &  1.6667 & 1.6667 \\
		(5,6) &  1.6667 & 1.6667 \\
		\hline
	\end{tabular}
	\caption{Values of $r_{n,n+1}$ for the Sierpinski triangle, given by \eqref{definicionrn}.}
	\label{tablarn}
\end{table}

%$$r_{3,4}= 1.342, \quad r_{4,5}=1.438, \quad r_{5,6}=1.513443, \quad r_{6,7}=1.598989.$$

\subsection{Renormalized finite elements methods- rFEMs}\label{MEF}

 Now we construct approximations for the weak form 
 \eqref{weakformJG}. That is, we propose approximations for the computations of renormalized energy bilinear forms by rescaling standard approximations as introduced in Section \ref{procedure}.
 
We defined $P^n=\{  u:V_n\rightarrow \mathbb{R} \}$.  We then project the weak formulation into the 
space  $P^n$. In the Galerkin formulation we seek to find  $u\in P^n$ such that
\begin{equation}
\begin{cases}
\mathcal{E}(u,v)=\displaystyle \int g v d\mu \mbox{ for all } v\in P^n\cap H_0^1,\\
u(x)=h(x), \quad x\in V_0.
\end{cases}
\label{galer}
\end{equation}
Both  $u$ and $v$ are defined on $V_n$, then this formulation is equivalent to 
\begin{equation}
\begin{cases}
\mathcal{E}_n(u,v)=\displaystyle \int g v d\mu \mbox{ for all } v\in P^n\cap H_0^1,\\
u(x)=h(x), \quad x\in V_0.
\end{cases}
\label{galer2}
\end{equation}
We propose several classical finite element methods procedures in one and two dimensions to approximate the renormalized energy. Note that we can use the energy  $\mathcal{E}_n$ rather than the limit $\mathcal{E}$. The renormalization constant is then computed by a procedure similar to the one explained previously.

\subsubsection{Integrals along the edges}\label{seccontornos}
Recall that  $\mathcal{E}_n(u,v)=-\int \Delta_n u\cdot v d\mu_n^{1}$ when $u,v:V_n\rightarrow \mathbb{R}$. Introduce the bilinear form, 
\begin{equation}
\mathcal{E}_n^{(1)}(u,v) = \sum_{x\sim^n y}\int_x^y u'v'd\nu_n^1,
\end{equation}
where $u'$  denotes the one dimensional derivative along $[x,y]$ of the linear interpolation of the vertex values   $u(x)$ and $u(y)$. 
The measure $\nu_n^1$ is defined as the length measure along the edges of $W_n$ (rescaled to obtain total length 1). Note that for each $n$ there are
total $3^{n+1}$ edges (3 for each cell). Each edge of $V_n$ has  length 
$2^{-n}$. Given a total length of $3 (3/2)^n$ before length rescaling.\\

We consider the following discrete problem,

\begin{equation}
\begin{cases}
2\mathcal{E}^{(1)}_n(u,v)=\displaystyle \int g v d\nu_n^1 \mbox{ for all  } v\in P^n\cap H_0^1.\\
u(x)=h(x), \quad x\in V_0.
\end{cases}
\label{galer2longitud}
\end{equation}
As before we can write  $u=u_I+u_G$ where $u_G(x)=g(x)$, $x\in V_0$; $u_G(x)=0$,  $x\in V_n\setminus V_0$. Analogously $u_I\in P^n$ tal such that $u_I(x)=0$ for 
$x\in V_0$. We then have, 
$$\mathcal{E}_n^{(1)}(u_I,v)d\nu_n^1=\int gd\nu_n^1-\mathcal{E}_n^{(1)}( 
u_G,v).$$
This is equivalent  to the linear system, 
$$A^{(1)}_nu=b^{(1)}_n.$$
Let  $V_n\setminus V_0=\{x_1,x_2,...,x_p\}$ the set of interior vertices and we have, 
$$a_{ij}= \sum_{x\sim^n y}\int_x^y\varphi_{x_i}'\varphi_{x_j}'d\nu_n^1,$$
where $\varphi_{x_i}$ is the linear interpolation of the characteristic functions of $\{x_i\}$ in $V_n$.
We also have
$$b_i=\int g\cdot \varphi_j d\mu_n -\int u_G'\varphi_j'd\nu_n^1.$$

It is easy to see that 
$$\int g \varphi_j d\nu_n^1=\sum_{a\sim_n b}\int_0^{2^{-n}} \left(2^n(g(b)-g(a))x-g(a)\right)\cdot \varphi_j(x)dx.$$
Recall that, on the left hand side above we use the piecewise liner interpolation of the nodal value of $g$. The previous formulation have to be renomarlized to \begin{equation}
\begin{cases}
(r^{(1)})^n\mathcal{E}^{(1)}_n(u,v)=\displaystyle \int g v d\nu_n^{(1)} \mbox{ for all } v\in P^n\cap H_0^1.\\
u(x)=h(x), \quad x\in V_0.
\end{cases}
\label{galer2longitudrenorm}
\end{equation}

We note the following relation between the renormalized measure and the measure induced by the length measure, for any $n$, we have 
\begin{equation}\label{relmeasures}
\int fd\mu_n\approx  2^{-n}\int fd\nu_n^{(1)}, 
\end{equation}
that follows by computing the length integrals using the trapezoidal rule.

\begin{rem}
The renormalization constant $r^{(1)}$ can be approximated using the procedure described 
in Section \ref{procedure} by solving consecutive refinement level approximations with the constant function 1 as the right hand side and Dirichlet boundary conditions. See \eqref{definicionrn}.
In Table 
\ref{tablarn2} we show the results of computing the renormalization constant.\\
\end{rem}

In this case the renormalization constant can also be computed analytically. Integration by parts and the fact that 
we use piecewise linear interpolation ($u''=0$ inside edges)
yields 
$$\mathcal{E}_n^{(1)}(u,v)=\sum_{x\sim_n y}  u'(t)v(t)|_x^y-\int_x^y u''(t)v(t)d\nu_n^1=
\sum_{x\sim_n y}  u'(y)v(y)-u'(x)v(x).$$
Having into account that the length of the  edges of the $n$ approximation of $K$ is  $1/2^n$, 
we get, 
$$\mathcal{E}_n^{(1)}(u,v)= \sum_{x\sim_n y} \frac{u(y)-u(x)}{1/2^n}v(y)-\frac{u(y)-u(x)}{1/2^n}v(x)= 2^n\sum_{x\sim_n y} (u(y)-u(x))(v(y)-v(x)).$$
Therefore, 
$$ A^{(1)}_n=2^n\cdot \widetilde{\Delta}_n 
\quad \mbox{ and } \quad  
\Delta_n= \left(\frac{5}{2}\right)^n  A^{(1)}_n.$$
Due to (\ref{relmeasures}) and 
\eqref{rforEn} we see that renormalization constant for the family $A^{(1)}_n$ is $r^{(1)}=\frac{5}{4}.$ In Table 
\ref{tablarn2} we show the results of computing the renormalization constant using the procedure explained in Section \ref{procedure}. This a  numerical verification that for the case of the Sierpinski triangle the computation agrees with the exact value of the renormalization constant just derived.

\begin{table}
	\centering
	\begin{tabular}{ | c | c | c |}
		\hline
		$(n,n+1)$ & $\max_{x}r_{n,n+1}$ &
		$\mbox{mean}_xr_{n,n+1}$\\
		\hline
		(4,5) &  1.2500 & 1.2500 \\
		(5,6) &  1.2500 & 1.2500 \\
		(6,7) &  1.2500 & 1.2500 \\
		\hline
	\end{tabular}
	\caption{Values of $r_{n,n+1}$ for the Sierpiski triangle using integrals along Edges.}
	\label{tablarn2}
\end{table}

%Igualando con la forma de la energ\'ia para el grafo, se obtiene que 
%$$k2^{n}=\left(\frac{5}{3}\right)^n.$$
%y de esta forma se llega al siguiente resultado
%$$k=\left(\frac{5}{6}\right)^n.$$
%Que esta va a ser la asociaci\'on entre la forma bilineal unidimensional, con la forma del Laplaciano asociada a un grafo.

%As\'i comparando las soluciones por el MEF a traves del contorno se obtiene que:
%
%$$a_c(u,v)=\left(\frac{5}{6}\right)^n \_n(u,v).$$ 
%
%Teniendo en cuenta el estudio anterior, la cual nos asocia la parte discreta con la parte continua sobre el contorno de los tri\'angulos, se tiene la siguiente equivalencia:
%
%$$\frac{5}{6}[a_c(f_i,f_j)]=[e_n(f_i,f_j)].$$

\subsubsection{Area integrals}

Introduce the bilinear form, 
\begin{equation}
\mathcal{E}_n^{(2)}(u,v) = \sum_{\tau\in K_n }\int_\tau \nabla u\nabla vd\nu_n^{(2)},
\end{equation}
where $\nabla u$  denotes the two-dimensional gradient of the two-dimensional linear interpolation of the nodal value of $u$ in the 
triangle $\tau$. 
The measure $\nu_n^{(2)}$  is the area  measure restricted to  $K_n$ and normalized such that 
the total area of (all the triangles of) $V_n$ is one. Note that, for each $n$, there is $3^n$ each of them of area $\frac{\sqrt{3}}{4}\frac{1}{2^{2n}}$ for a total area of $\frac{\sqrt{3}}{4}(\frac{3}{4})^n$ before rescaling.\\

We formulate the following discrete problem, 
\begin{equation}
\begin{cases}
\mathcal{E}^{(2)}_n(u,v)=\displaystyle \int g v d\nu_n^{(2)} \mbox{ para todo } v\in P^n\cap H_0^1.\\
u(x)=h(x), \quad x\in V_0.
\end{cases}
\label{galer2longitud}
\end{equation}

This time the previous formulation is equivalent to the linear system, 
$$A^{(2)}_n u=b^{(2)}_n, $$
where
$$a_{i,j}=\int \nabla\varphi_i\nabla\varphi_jd\nu_n^{(2)},$$
and 
$$b_i=\int f\varphi_id\nu_n^{(2)}.$$

As before, a renormalization is needed, that is,
\begin{equation}
\begin{cases}
(r^{(2)})^n\mathcal{E}^{(2)}_n(u,v)=\displaystyle \int g v d\nu_n^{(2)} \mbox{ para todo } v\in P^n\cap H_0^1.\\
u(x)=h(x), \quad x\in V_0.
\end{cases}
\label{galer2areasrenorm}
\end{equation}

\begin{rem}
The renormalization constant $r^{(2}$ can be approximated using the procedure described 
in Section \ref{procedure} by solving consecutive refinement level approximations with the constant function 1 as the right hand side and Dirichlet boundary conditions. See \eqref{definicionrn}.
In Table 
\ref{tablarn3} we show the results of computing the renormalization constant.\\
\end{rem}

In order to verify our computations we compute 
the renormalization constant analytically. This is possible in this case. Recall that, 
$$\mathcal{E}^{(2)}_n(f,g)= \sum_{\tau\in K_n }\int_\tau \frac{\partial f}{\partial x}\frac{\partial g}{\partial x}+\frac{\partial f}{\partial y}\frac{\partial g}{\partial y}\quad d\nu_n^2.$$
We use standard finite element analysis. Introduce the reference basis functions 
$$\hat{P}_1(\hat x,\hat y)=1-\hat{x}-\hat{y}, \quad \hat{P}_2(\hat x,\hat y)=\hat x, \quad \hat{P}_3(\hat x,\hat y)=\hat y.$$
defined in the reference triangle $\hat{\tau}$ with vertices $(0,0)$, $(1,0)$ and $(0,1)$. This reference triangle can be mapped into the triangles of $K_n$ by an affine mapping in two dimensions. If we consider the triangle $\tau$ of $K_n$ with vertices
$\vec{x}_1=(x_1,y_1), \vec{x}_2=(x_2,y_2)$ and $\vec{x}_3=(x_3,y_3)$, this mapping is given by 
$F_\tau:  \hat{\tau}\to \tau$ defined by
$$F(\vec{x})=V\vec{x} +\vec{x}_1; \quad V=\left(\begin{array}{cc}
x_2-x_1 & y_2-y_1\\\\
x_3-x_1 & y_3-y_1
\end{array}\right).$$
Define $P_i(\vec{x})=\hat{P}_i(F^{-1}(\vec{x}))$, $i=1,2,3$. Any linear function on $\tau$ is a linear combination of the basis functions $P_1,P_2,P3$, in particular if $u$ is a linear function on 
$\tau$ we have $u(\psi)=u(x)P_1(\psi)+u(y)P_2(\psi)+u(z)P_3(\psi)$.
From the definition of $P_i$ is easy to see that 
$$V^{-1}=\left(\begin{array}{cc}
\displaystyle \frac{\partial P_2}{\partial x} &  \displaystyle  \frac{\partial P_3}{\partial x}\\ \\
\displaystyle  \frac{\partial P_2}{\partial y} & \displaystyle  \frac{\partial P_3}{\partial y}
\end{array}\right)=\frac{1}{\det(V)}\left(\begin{array}{cc}
y_3-y_1 & y_1-y_2\\
x_1-x_3 & x_2-x_1\\
\end{array}\right).$$
Moreover, we also have, 
$$\frac{\partial P_1}{\partial x}=-\frac{\partial P_2}{\partial x}-\frac{\partial P_3}{\partial x}, \qquad \frac{\partial P_1}{\partial y}=-\frac{\partial P_2}{\partial y}-\frac{\partial P_3}{\partial y}.$$
We also recall that $\det(V)=\frac{\sqrt{3}l^2}{2}$
where  $l=2^{-n}$ is the diameter of the triangle. 
We can then compute, 
$$\begin{array}{rcl}
a_A(P_1,P_1) & = & \displaystyle \int_\tau {\left(\frac{\partial P_1}{\partial x}\right)^2+\left(\frac{\partial P_1}{\partial y}\right)^2 d\nu^{(2)}}\\
\\
& = & \displaystyle \frac{4}{3l^4} \int_\tau (y_3-y_2)^2+(x_3-x_2)^2 d\nu^{(2)}  \\
\\
& = &  \displaystyle \frac{4}{3l^2} \int_\tau 1 d\nu^{(2)} \quad  \mbox{since}\quad (y_3-y_2)^2+(x_3-x_2)^2=l^2,\\
& =&  \displaystyle \frac{\sqrt{3}}{3}. \end{array}$$
Analogously we have,
$$\begin{array}{rcl}
a_A(P_1,P_2) & = & \displaystyle\int_\tau{\frac{\partial P_1}{\partial x}\frac{\partial P_2}{\partial x}+\frac{\partial P_1}{\partial y}\frac{\partial P_2}{\partial y} d\nu^{(2)}}\\
& = & \displaystyle \frac{4}{3l^4} \int_\tau -l^2 + \langle  (x_2-x_1,y_2-y_1), (x_3-x_1,y_3-y_1) \rangle d\nu^{(2)}  \\
\\
& = & \displaystyle \frac{4}{3l^2} \displaystyle\int_\tau l^2(1-\cos(60^\circ) )d\nu^{(2)} 
=  \displaystyle - \frac{\sqrt{3}}{6}.
\end{array}$$
Having into account that each interior node belongs only to two-triangles and that two distinct nodes share an edge in at most one triangle, we conclude that the assembled global matrix is given by   
\begin{equation}
A_n^{2}=[a_{ij}]=\begin{cases}
\frac{\sqrt{3}}{3} \quad \mbox{ if } i=j.\\
-\frac{\sqrt{3}}{6} \quad \mbox{ if } i\neq j,
\end{cases}
\end{equation}
where $i$ and $j$ corresponds to the index of interior nodes. We see that, 
\begin{equation}
A_n^{(2)}=[a_{ij}]=\frac{\sqrt{3}}{6}\cdot\begin{cases}
p_{ii} \quad \mbox{ if } i=j.\\
-1 \quad \mbox{ if } i\sim_n j.
\end{cases}, \quad \mbox{ and then }
\quad A_n^{(2)}=\frac{\sqrt{3}}{6}\cdot\widetilde{\Delta}_n.
\label{areasformula}\end{equation}
Recall that $\widetilde{\Delta}_n$ was defined as the graph laplacian of the 
 $n$ approximation of the Sierpinski triangle.\\  
 
We note the following relation between the renormalized measure and the measure induced by the area measure, for any $n$, we have 
\begin{equation}\label{relmeasures2}
\int fd\mu_n\approx \frac{4}{\sqrt{3}} 4^{-n}\int fd\nu_n^{(2)}, 
\end{equation}
that follows by computing the length integrals using the trapezoidal rule.

From (\ref{areasformula}) and (\ref{relmeasures2}), we then have that the exact value of the renormalization constant is the same as the one for the construction based on length measures. This result verifies the numerical computations obtained in Table \ref{tablarn3}.

\begin{table}
	\centering
	\begin{tabular}{ | c | c | c |}
		\hline
		$(n,n+1)$ & $\max_{x}r_{n,n+1}$ &
		$\mbox{mean}_xr_{n,n+1}$\\
		\hline
		(4,5) &  1.2500 & 1.2500 \\
		(5,6) &  1.2500 & 1.2500 \\
		(6,7) &  1.2500 & 1.2500 \\
		\hline
	\end{tabular}
	\caption{Values of $r_{n,n+1}$ for the Sierpiski triangle using  area measures.}
	\label{tablarn3}
\end{table}
%cuando $i=j$, que ser\'ia en la diagonal de la matriz laplaciana del grafo es 2 en los v\'ertices del tri\'angulo $\tau_0$, haciendo c\'alculos sobre eso v\'ertices, se tiene
%$$\begin{array}{rcl}
%\frac{2\sqrt{3}}{3}\frac{1}{2^{2n}}h & = & \left(\frac{5}{3}\right)^n (2) \\
%\\
%h & = & \left(\frac{20}{3}\right)^n \frac{3}{\sqrt{3}}
%\end{array}$$
%esta ser\'ia la constante que asociar\'ia $a_A$ con $a_e$

%
\subsection {Illustrations of the numerical methods}
This section shows the numerical solution of the Laplace equation posed in some fractal sets. In particular, we consider the Sierpinski triangle, the Kosh curve, and two Hata trees. As we discussed before, we can approximate the  solution by\begin{itemize}
\item Renormalized Finite Difference (rFD): 
\begin{itemize}
\item Pre-processing: We solve a model problem with a given Dirichlet condition and $g(x)=1$ as the forcing term in several graph approximations of the fractal set to compute an approximation of the renormalization constant.
\item Online step: we solve for the actual forcing term with the renormalized graph laplacian using the approximation of the renormalization constant computed in the pre-processing step.
\end{itemize}
\item Renormalized Finite Element Method with line integrals (rFEM1D): 
\begin{itemize}
\item Pre-processing: approximation of the renormalization constant as before.
\item Online step: solution with actual right-hand side
\end{itemize}
\item Renormalized Finite Element Method with area integrals (rFEM2D): 
\begin{itemize}
\item Pre-proccesing: approximation of the renormalization constant as before.
\item Online step: solution with actual right hand side
\end{itemize}

\end{itemize}
\subsubsection{The Sierpinski triangle}
\begin{figure}[H]
	
	\centering 
	\includegraphics[scale=0.5]{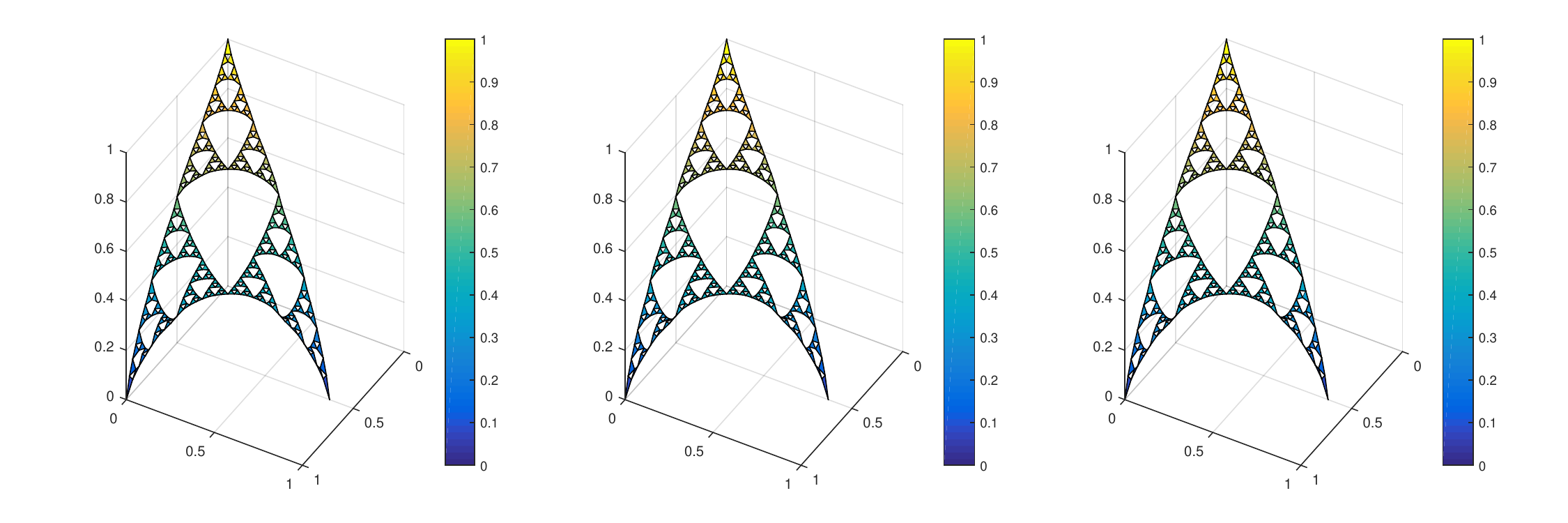}
	\caption{Approximated solution with  rFD (left), rFEM1D (center) and rFEM2D (right). Here $g(x,y)=\sin(x+y)$, $u(a_0)=1$,  $u(a_1)=u(a_2)=0$ and we consider $V_5$ the fifth level approximation of the Sierpinski triangle.}
	\label{solsierc1}
\end{figure}

For problems posed on the Sierpinski triangle recall that 
$V_0=\{a_0,a_1,a_2\}$. We want to approximate the solution of
$$\left\{ \begin{array}{lr}
-\Delta u(x,y)=f(x,y),\\
\\
u(a_0)=1, \quad u(a_1)=0, \quad u(a_2)=0.\\
\end{array}\right.
$$
In Figure \ref{solsierc1} we illustrate some results.

\subsubsection{The Kosh curve}
For problems posed on the Kosh curve, recall that $V_0=\{(0,0),(1,0)\}$. We want to approximate the solution of

$$\left\{ \begin{array}{lr}
-\Delta u(x,y)=g(x,y),\\
\\
u(a_0)=1, \quad u(a_1)=0.\\
\end{array}\right.
$$
We use the method introduced in section  \ref{seccontornos}; if we associate the Laplacian with the energy defined by integrations over the edges, we have 
\begin{equation}
\begin{cases}
(r^{(1)})^n\mathcal{E}^{(1)}_n(u,v)=\displaystyle \int g v d\mu_n \mbox{ for all } v\in P^n\cap H_0^1.\\
u(x)=h(x), \quad x\in V_0,
\end{cases}
\label{igualdadconmedida}
\end{equation}
where $\mu_n$ is the length measure restricted to edges and   $r^{(1)}=\frac{16}{9}$ was computed in  Table  (\ref{tablarnK}) using the procedure 
described in Section \ref{procedure}.
\begin{figure}[H]
	\centering 
	\includegraphics[scale=0.4]{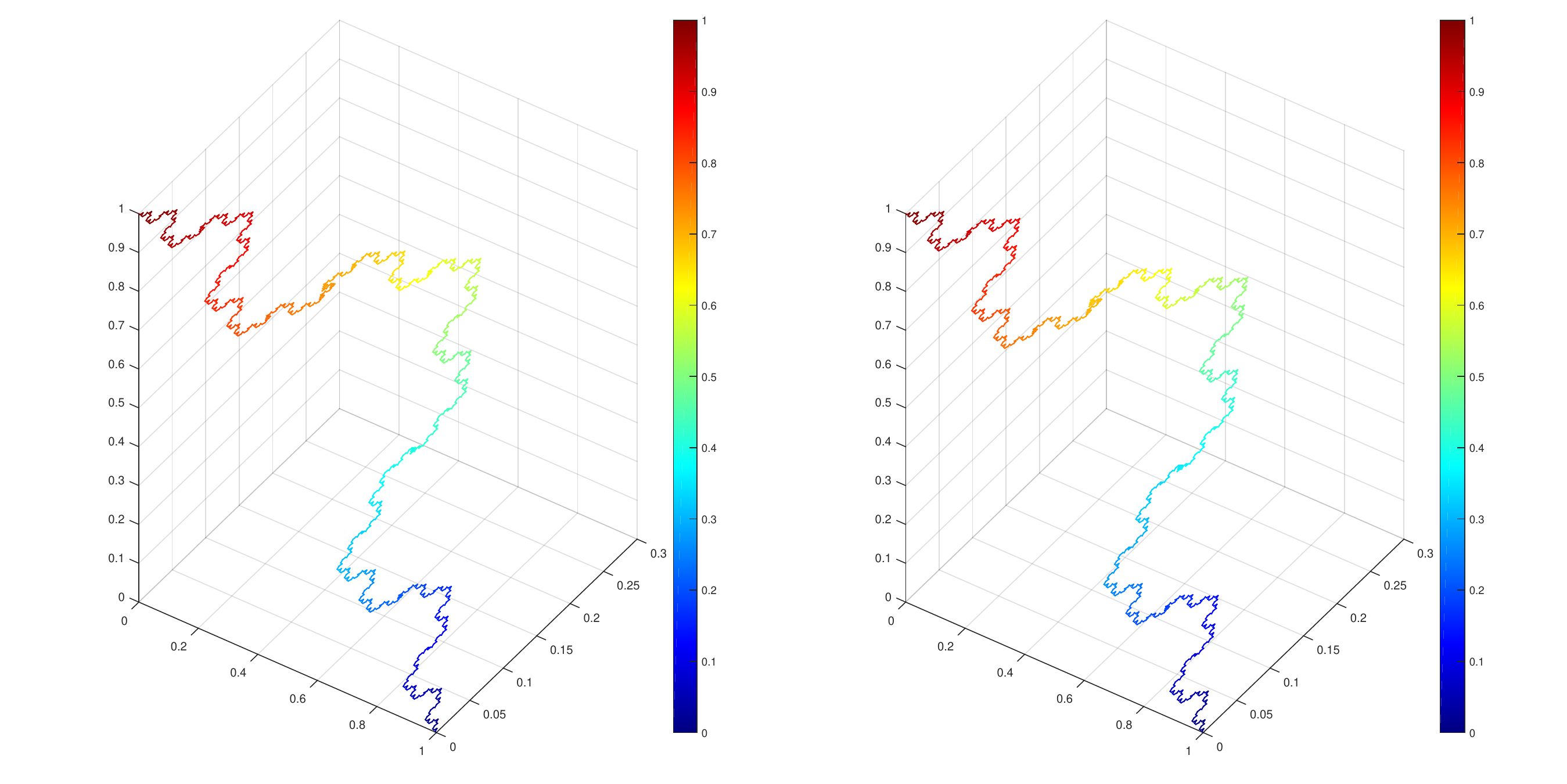}
	\caption{\small Approximated solution with  rFD (left) and rFEM1D (right). Here $g(x,y)=\sin(x+y)$, $u(a_0)=1$,  $u(a_1)=0$ and we consider $V_5$ the fifth level approximation of the Kosh curve.}
	\label{solkosh}
\end{figure}

\begin{table}
	\centering
	\begin{tabular}{ | c | c | c |}
		\hline
		$(n,n+1)$ & $\max_{x}r_{n,n+1}$ &
		$\mbox{mean}_xr_{n,n+1}$\\
		\hline
		(3,4) &  1.7778 & 1.7778 \\
		(4,5) &  1.7778 & 1.7778 \\
		(5,6) &  1.7778 & 1.7778 \\
		\hline
	\end{tabular}
	\caption{Values of $r_{n,n+1}$ that the Kosh curve.}
	\label{tablarnK}
\end{table}

\begin{figure}[H]
	
	\centering 
	\includegraphics[scale=0.4]{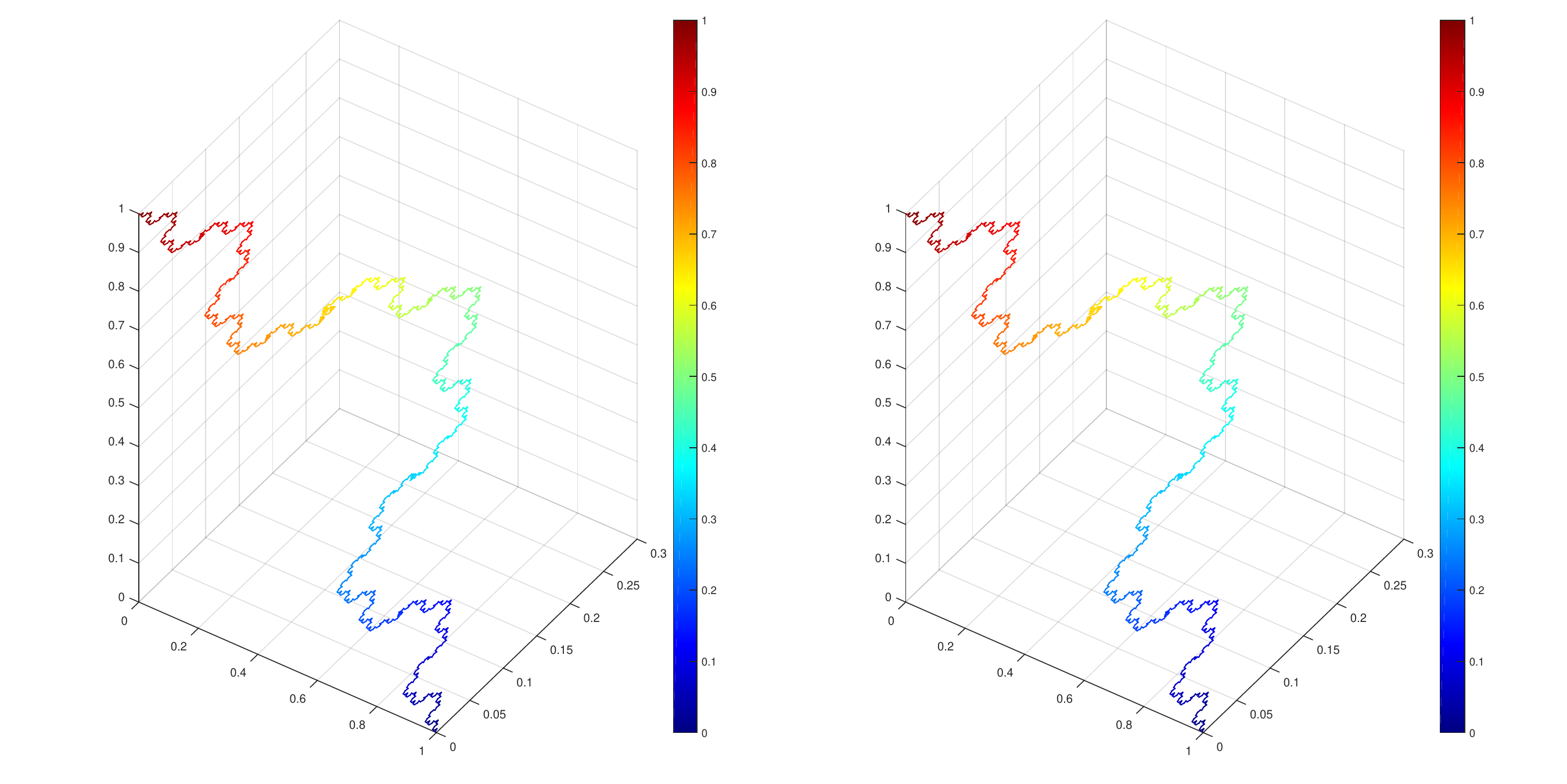}
	\caption{\small Approximated solution with  rFD (left) and rFEM1D (right). Here $g(x,y)=0$, $u(a_0)=1$,  $u(a_1)=0$ and we consider $V_3$ the third level approximation of the Kosh curve.}
	\label{solkosha}
\end{figure}

\subsubsection{The Hata tree}
In the find method, the solution of the following equation is the same as that studied before, where $V_0=\{(0,0),(1,0)\}$.

$$\left\{ \begin{array}{lr}
-\Delta u(x,y)=g(x,y).\\
\\
u(a_0)=1, \quad u(a_1)=0,\\
\end{array}\right.
$$

if we associate the Laplacian with the energy on edge, we have

\begin{equation}
\begin{cases}
(r^{(1)})^n\mathcal{E}^{(1)}_n(u,v)=\displaystyle \int g v d\mu_n \mbox{ for all } v\in P^n\cap H_0^1.\\
u(x)=h(x), \quad x\in V_0.
\end{cases}
\label{igualdadconmedida}
\end{equation}

Where $\mu_n$ is the self-similar measure that the Kosh curve and the renormalization constant is given for  $r^{(1)}=\frac{5}{3}$, see Table  \ref{tablarnK}.

\begin{table}
	\centering
	\begin{tabular}{ | c | c | c |}
		\hline
		$(n,n+1)$ & $\max_{x}r_{n,n+1}$ &
		$\mbox{mean}_xr_{n,n+1}$\\
		\hline
		(3,4) &  1.6667 & 1.6667 \\
		(4,5) &  1.6667 & 1.6667 \\
		(5,6) &  1.6667 & 1.6667 \\
		\hline
	\end{tabular}
	\caption{Values of $r_{n,n+1}$ for the Hata tree.}
	\label{tablarnK}
\end{table}

\begin{figure}[H]
	\centering 
	\includegraphics[scale=0.4]{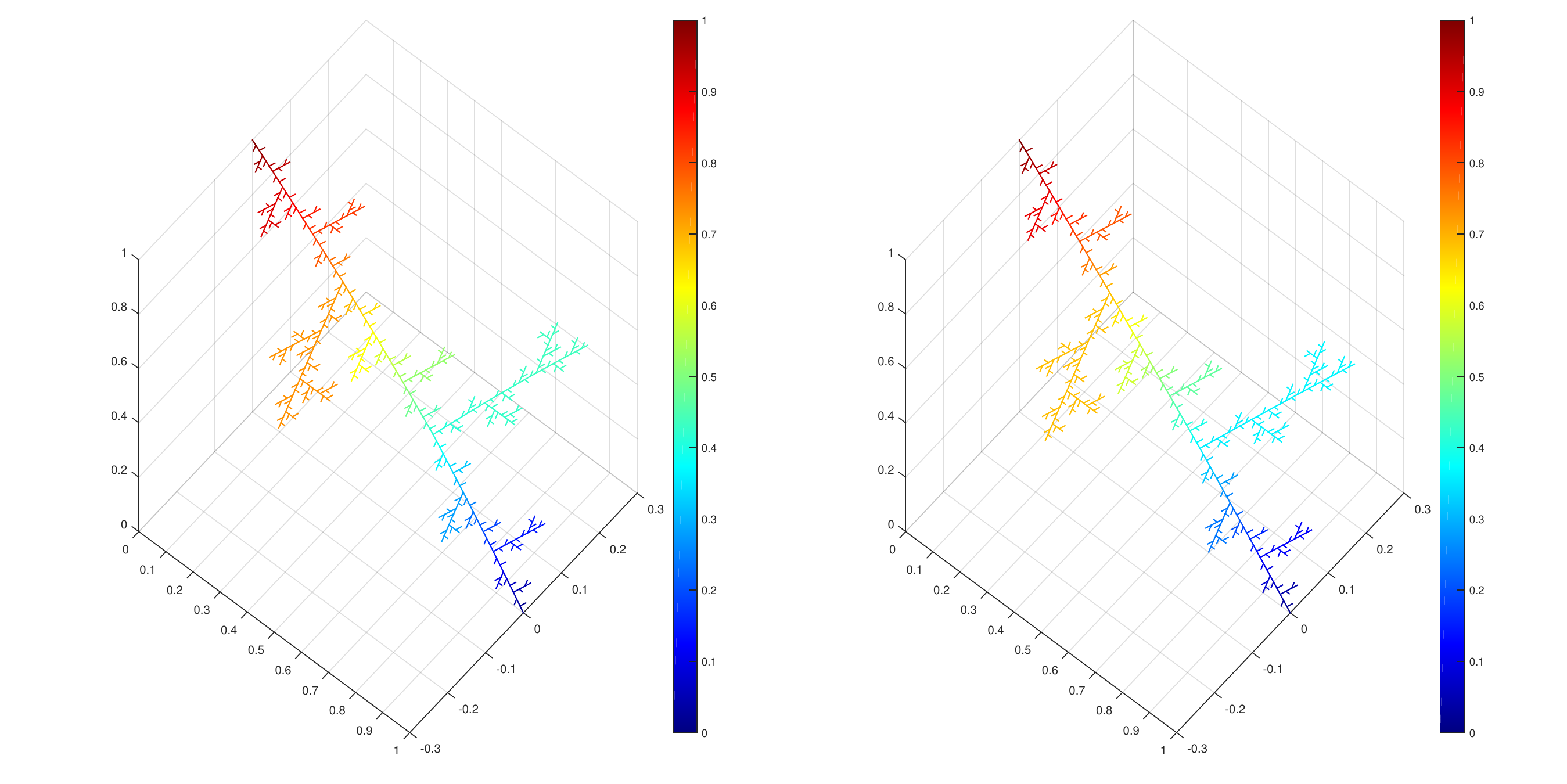}
	\caption{\small Approximated solution with  rFD (left) and rFEM1D (right). Here $g(x,y)=\sin(x+y)$, $u(a_0)=1$,  $u(a_1)=0$ and we consider $V_3$ the third level approximation of the Hata tree.}
	\label{solhata}
\end{figure}

\begin{figure}[H]
	
	\centering 
	\includegraphics[scale=0.4]{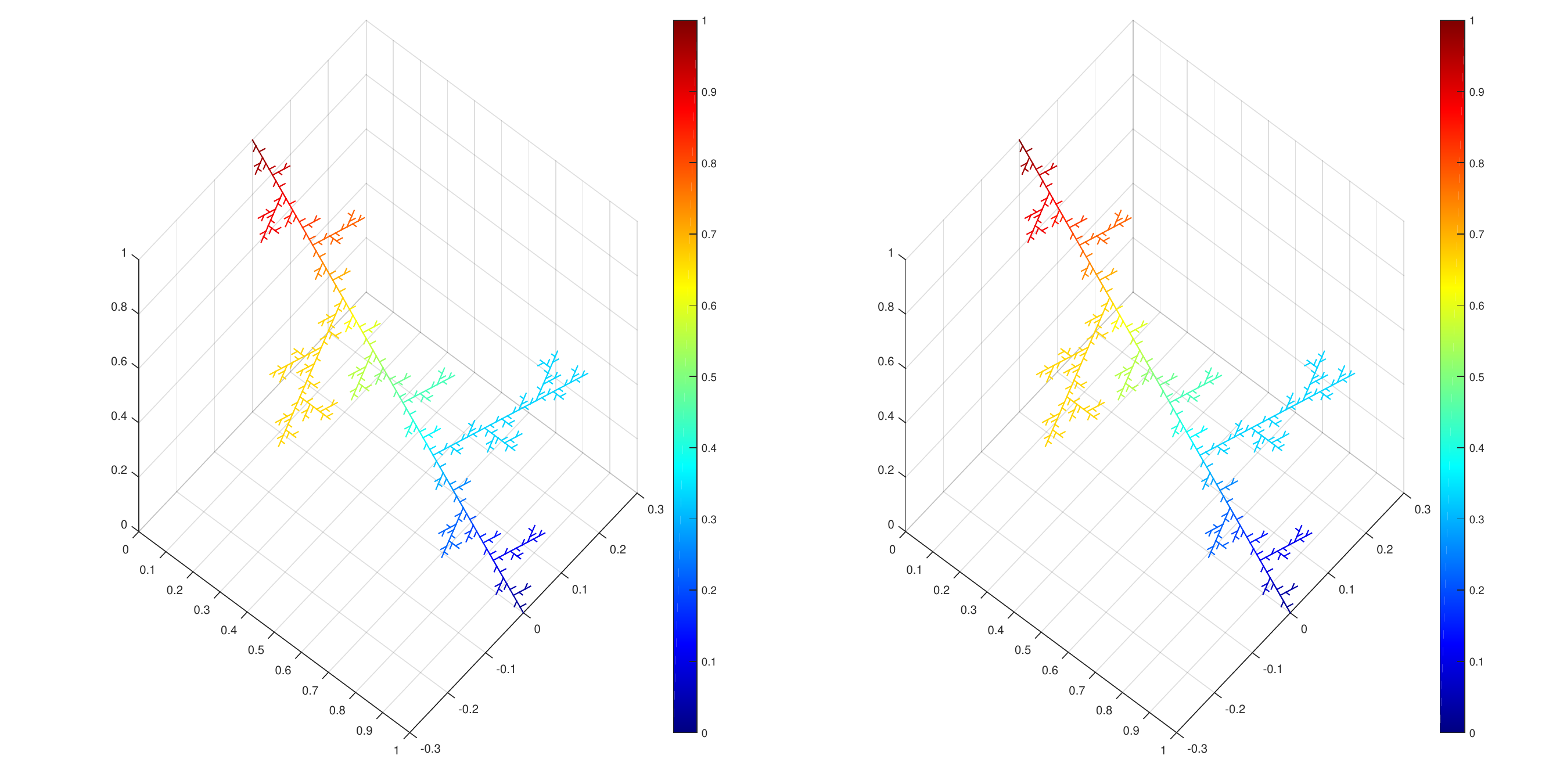}
	\caption{\small Approximated solution with  rFD (left) and rFEM1D (right). Here $g(x,y)=0$, $u(a_0)=1$,  $u(a_1)=0$ and we consider $V_3$ the third level approximation of the Hata tree.}
	\label{solhataa}
\end{figure}

\subsubsection{Hata tree in the space}
\begin{figure}[H]
	\centering 
	\includegraphics[scale=1.1]{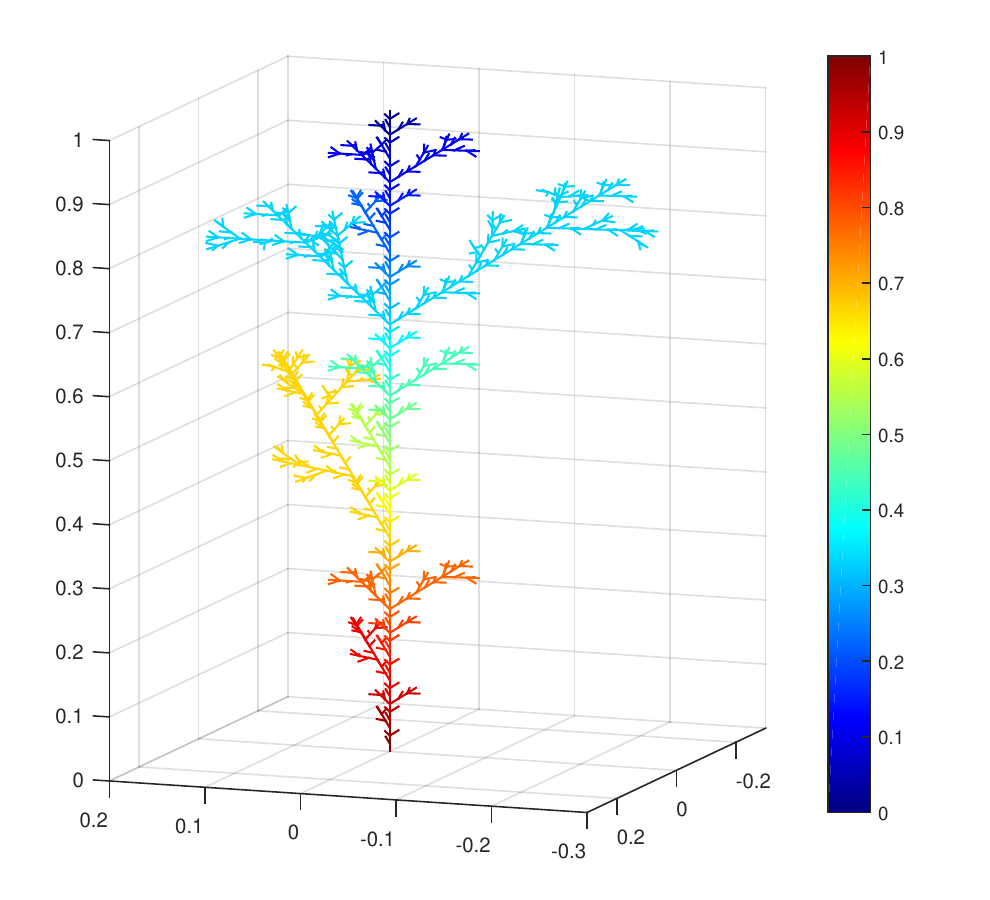}
	\caption{\small Approximated solution with  rFD. Here $g(x,y)=0$, $u(a_0)=1$,  $u(a_1)=0$ and we consider $V_3$ the third level approximation of the Hata tree in the space.}
	\label{solhata3d}
\end{figure}

For this fractal we use only rFD method (\ref{diferenciafinitas}), In Figure \ref{solhata3d} we present a computed solution.

\section{Conclusions}\label{sec:conclusions}

This paper designed a numerical procedure to approximate solutions to diffusion problems on self-similar fractal sets. We start with a discrete approximation of the fractal and the derivatives in standard non-renormalized formulations. We can then precompute the renormalization constant needed to approximate the actual differential operators on the fractal set. In particular, we present examples with the Sierpinski 
triangle using standard graph weights and adjacency matrices (Finite Difference method) or using week forms with length or area measures (Finite Element method). In the Finite Element method with length measure, the derivatives in the weak forms are classical derivatives along the edges with integration concerning the length measure. In the finite element method with area measure, we use partial derivatives with integration in two dimensions on triangles of the approximation of the Sierpinski triangle. We also present additional illustrations with the Kosh curve and the Hata tree.

It is also important to mention that the implementation of finite elements is simple and does not have significant changes to the finite element method for differential equations in open domains. Also, the renormalization constant does not need to be known a priori. We can use finite element codes that work on triangulations in general, and only the `` triangulation'' or graph that approximates the fractal must be used as input for these codes. The renormalization constant can be precomputed as proposed in this paper. We observe that diffusion processes on these fractal sets can be approximated by classical diffusion processes (involving classical derivatives) on fractal approximations. These must be rescaled by the scale parameter that can be precomputed. The authors will explore this idea and the related numerical analysis in future works. \\

\section*{Acknowledgements}
The authors thank Professor Milton Jara for introducing us to the topic of diffusion on fractals.

\end{document}